# ISOSCATTERING SCHOTTKY MANIFOLDS

ROBERT BROOKS, RUTH GORNET, AND PETER PERRY

ABSTRACT. We exhibit pairs of infinite-volume, hyperbolic three-manifolds that have the same scattering poles and conformally equivalent boundaries, but which are not isometric. The examples are constructed using Schottky groups and the Sunada construction.

## 1. INTRODUCTION

The purpose of this paper is to explore the spectral geometry of the scattering operator by exhibiting examples of infinite volume hyperbolic three-manifolds that are 'isoscattering' in a sense we will make precise, but have distinct geometries. To do so we will work with convex co-compact hyperbolic manifolds associated to Schottky groups of hyperbolic isometries.

This class of manifolds, described in greater detail in what follows, may be thought of as interiors $X$ of handlebodies $\bar{X}$ equipped with a metric that puts the boundary at metric infinity. The interior carries a hyperbolic structure, and the boundary, $\partial \bar{X}$, carries an induced conformal structure as a Riemann surface.

Associated to the hyperbolic structure is a Laplace-Beltrami operator with at most finitely many $L^2$ eigenvalues and continuous spectrum of infinite multiplicity [22, 23]. The scattering operator, $S_X(s)$, characterizes asymptotic behavior of generalized eigenfunctions associated to the continuous spectrum. It is a pseudo-differential operator on $C^\infty(\partial \bar{X})$. The relation $S_X(s)f_- = f_+$ may be thought of informally as follows: a wave energy $s(2-s)$ and initial amplitude $f_-$ 'at infinity' propagates into the interior of $X$, scatters, and re-emerges as a wave with final amplitude $f_+$ 'at infinity.' The scattering operator depends meromorphically on the spectral parameter $s$.

Roughly and informally, the scattering poles of $S_X(s)$ arise because some of the energy of a scattered wave can become temporarily 'trapped' near closed geodesics of $X$. The poles of $S_X(s)$ are known to be among the zeros of the zeta function for geodesic flow on $X$ [37], and it is a 'folk theorem' that a Selberg-type trace formula relates the closed geodesics and scattering poles of $X$ much as the celebrated Selberg trace formula for compact Riemann surfaces [42], and its generalization to compact manifolds by Duistermaat-Guillemin [13] relate the closed geodesics and eigenvalue spectrum of the Laplacian (see Colin de Verdiére [12] for similar results using the heat kernel, see Buser [11] for further discussion of the compact Riemann surface case). Moreover, in simple models of infinite volume hyperbolic manifolds, it is easy to see how a closed geodesic gives rise to an infinite sequence of scattering

Date: October 3, 1998.

R. B. supported in part by the Israel Academy of Sciences, the Fund for the Promotion of Research at the Technion, and the New York Metropolitan Fund.

R. G. supported in part by NSF Grant DMS 97-53220.

P. P. supported in part by NSF Grant DMS 97-07051.

poles (e.g. [35], Appendix B). Thus it is natural to regard the scattering poles as analytic data analogous to the eigenvalues of the Laplacian on a compact surface and to investigate their geometric content (see [47] for a survey and [17, 18, 19] for more recent results).

We will call two Schottky manifolds $\bar{X}_1$ and $\bar{X}_2$ *isoscattering* if the scattering operators $S_{X_1}(s)$ and $S_{X_2}(s)$ have the same scattering poles counted with multiplicities, as defined in §2. We will call them *strongly isoscattering* if there is an invertible linear map $T_\vartheta : C^\infty(\partial \bar{X}_1) \rightarrow C^\infty(\partial \bar{X}_2)$, independent of $s$, that intertwines the scattering operators $S_{X_1}(s)$ and $S_{X_2}(s)$. As the terminology suggests, if two Schottky manifolds are strongly isoscattering, they are also isoscattering.

Our main result is:

**Theorem 5.1.** *There exist $\bar{X}_1$ and $\bar{X}_2$ such that*

(i) *$X_1$ is strongly isoscattering to $X_2$,*

(ii) *$\partial \bar{X}_1$ is conformally equivalent to $\partial \bar{X}_2$, and*

(iii) *$X_1$ is not isometric to $X_2$.*

The result is proved using Bérard's version of the Sunada construction [2] together with an explicit construction of Schottky manifolds using ideas from Riemann surface theory. The Sunada construction was previously used to construct isoscattering surfaces by Bérard [2], Guillopé-Zworski ([19], Remark 2.15), and Zelditch [46].

Another simple set of three-dimensional isoscattering manifolds may be obtained as follows. Sunada's construction can be used to obtain isospectral Riemann surfaces of genus 4 and higher (see [9] and references therein). Let $S_1$ and $S_2$ be two such surfaces. The manifold $X_i = (0,1) \times S_i$ can be given a hyperbolic warped product metric and are conformally compact with $\bar{X}_i = [0,1] \times S_i$. Moreover (e.g. [35]), the scattering poles of $X_i$ are computable entirely in terms of the eigenvalues and Euler characteristic of the $S_i$ equipped with the hyperbolic metric. Thus $X_1$ and $X_2$ are isoscattering. To our knowledge, ours is the first set of counterexamples involving three-dimensional hyperbolic manifolds of infinite volume other than these.

In what follows, we will first discuss spectral and scattering theory for a class of hyperbolic manifolds that includes Schottky manifolds (§2). Next, we recall some basic facts about Schottky manifolds that will be important in our construction (§3). In §4 we recall Bérard's transplantation version [2] of Sunada's construction and prove a transplantation result that is then used in §5 to construct the examples.

**Acknowledgements.** Robert Brooks gratefully acknowledges the hospitality of the University of Kentucky. Ruth Gornet acknowledges the hospitality of the University of Kentucky where her work was supported by Texas Tech University and the National Science Foundation through a POWRE Fellowship. Peter Perry gratefully acknowledges the hospitality of the Department of Mathematics, Technion-Israel Institute of Technology, where part of this work was carried out.

## 2. Scattering Theory

In this section we recall some basic facts about geometry and spectral theory on convex co-compact hyperbolic manifolds. Let $\Gamma$ be a discrete group of isometries of real hyperbolic space $\mathbb{H}^3$. The group $\Gamma$ is geometrically finite if it admits a finite-sided geodesic polyhedron as a fundamental domain. A necessary but not sufficient condition is that $\Gamma$ be finitely generated. The group $\Gamma$ also acts on the geometric boundary of $\mathbb{H}^3$ (visualized, for example, as the Riemann sphere in the ball model of hyperbolic space) by conformal transformations. This action partitions $S^2$ into the domain of discontinuity, $\Omega(\Gamma)$, the largest domain on which $\Gamma$ acts properly discontinuously, and $\Lambda(\Gamma)$, the limit set for $\Gamma$. The group $\Gamma$ is called convex co-compact if the hyperbolic convex hull of the limit set intersects any fundamental domain in a compact subset of $\mathbb{H}^3$; this condition excludes fundamental domains with cusps. If $\Gamma$ is also torsion-free–which will be true for the examples we construct–the orbit space $X = \Gamma \setminus \mathbb{H}^3$ is a smooth manifold with a natural hyperbolic structure.

The manifold $X$ is the interior of a compact manifold with boundary $\bar{X}$, called the Klein manifold, given by $\bar{X} = (\Gamma \setminus \mathbb{H}^3) \cup (\Gamma \setminus \Omega(\Gamma))$. Thus $\partial \bar{X} = \Gamma \setminus \Omega(\Gamma)$ has a natural conformal structure. Moreover, $X$ is a *conformally compact manifold* [27, 30]: that is, the hyperbolic metric $g$ on $X$ takes the form $\rho^{-2} h$ where $\rho \in C^\infty(\bar{X})$ is a defining function for $\bar{X}$ (i.e., a smooth positive function on $\bar{X}$ that vanishes to first order at $\partial \bar{X}$) and $h$ is the restriction to $X$ of a smooth nondegenerate metric for $\bar{X}$. Note that $\rho$ and $h$ are determined up to multiplication by a smooth, bounded, invertible function on $\bar{X}$. We denote the Laplace-Beltrami operator on $X$ by $\Delta_X$.

The operator $\Delta_X$ is a positive operator with at most finitely many eigenvalues in $[0,1)$ (see [22], Theorem 4.8), and if the Hausdorff dimension of the limit set, $\Lambda(\Gamma)$, is less than one, there are *no* $L^2$ eigenvalues. The spectrum of $\Delta_X$ in $[1, \infty)$ is absolutely continuous (see [23]). Thus, the resolvent operator

$$R(s) = (\Delta_X - s(2-s))^{-1}$$

is meromorphic in the half plane $\Re(s) > 1$ with at most finitely many poles due to the eigenvalues of $\Delta_X$. The remaining spectral information is contained in the continuous spectrum, corresponding to the line $\Re(s) = 1$ in the complex $s$-plane, and encoded in the scattering operator. We wish to elucidate what geometric information is contained (or not contained) in the scattering operator.

To define the scattering operator, we recall the following uniqueness theorem [4, 21, 31]:

**Theorem 2.1.** *Let $\Re(s) = 1$ with $s \neq 1$, and let $f_- \in C^\infty(\partial \bar{X})$. Then there is a unique solution $u \in C^\infty(X)$ of the eigenvalue equation $(\Delta_X - s(2-s))u = 0$ having the asymptotic behavior*

$$u \sim \rho^s f_+ + \rho^{2-s} f_- + O(\rho^2)$$

*as $\rho \downarrow 0$, where $f_+ \in C^\infty(\partial \bar{X})$.*

This uniqueness theorem is analogous to the uniqueness theorem for harmonic functions in the upper half-plane with given boundary data, and implies that there is a one-to-one map $\mathcal{P}_s : C^\infty(\partial \bar{X}) \to C^\infty(X)$, the Poisson map, taking $f_-$ to $u$. The 'Poisson kernel' for this problem is a generalized eigenfunction $E_X(m, b, s)$ of $\Delta_X$ with eigenvalue $s(2-s)$. The uniqueness also implies that $f_+$ is determined by $f_-$, so there is a linear mapping $S_X(s) : C^\infty(\partial \bar{X}) \to C^\infty(\partial \bar{X})$ defined by $S_X(s)f_- = f_+$. This map is the scattering operator for $\Delta_X$ on $X$, and it follows from this definition that $S_X(s)S_X(2-s) = I$. Note that the scattering operator depends on a choice of defining function for $X$. This dependence is trivial since, if $\rho$ and $\hat{\rho}$ are defining functions for $\bar{X}$, $\hat{\rho} = h\rho$ for a bounded smooth function that is strictly positive on $\bar{X}$, and the corresponding scattering operators are related by $\widehat{S}(s) = h^{-s}S(s)h^{n-s}$. In our construction, we will make a choice for $\rho$ that is natural to the problem.

The scattering operator $S_X(s)$ extends to a meromorphic family of elliptic pseudodifferential operators on $\mathbb{C}$ with $S_X(s)$ of order $(2\Re(s)-2)$ [1, 15, 21, 24, 30, 36]. As defined above, the scattering operator has trivial poles in $\Re(s) > 1$ inherited from the scattering operator on the covering space $\mathbb{H}^3$. For this reason, it is convenient to work instead with the renormalized operator

$$\mathcal{S}_X(s) = 2^{2s-2}\frac{\Gamma(s-1)}{\Gamma(1-s)}S_X(s),$$

where $\Gamma(s)$ is the Gamma function. Thus $\mathcal{S}_X(s)$ has at most finitely many poles in the half-plane $\Re(s) > 1$. The poles of $\mathcal{S}_X(s)$ are among the zeros of Selberg's zeta function for the group $\Gamma$ [37] and thus carry geometric information. Moreover, it follows from the connection between scattering poles and zeros of the zeta function established in [36] and the lower bound on zeros of the zeta function established in [43], Example 6, p. 856 and Remark 3, p. 851, that the scattering operator has infinitely many poles in the half-plane $\Re(s) < 1$.

In line with the question "Can one hear the shape of a drum?" (i.e., do the eigenvalues of a compact manifold determine the manifold up to isometry; see for example the survey [3]) we would like to know to what extent $\mathcal{S}_X(s)$ determines $(X, g)$. It is known, for instance that if $\bar{X}_2$ is a quasi-conformal deformation of $\bar{X}_1$, and the difference of scattering operators is trace-class, then $X_1$ and $X_2$ are isometric [40]. Moreover, it is known that the scattering operator determines the conformal structure on $\partial\bar{X}$ through its principal symbol [5, 21]. Finally, it is shown in [5, 6] that if $X_1$ and $X_2$ are quasi-isometric then the norm of the difference of scattering operators measures the $K$-quasi isometry of $X_1$ and $X_2$.

Here we focus on the scattering poles as geometric data analogous to eigenvalues of the Laplacian on a compact manifold. We define the *multiplicity* of a scattering pole $s_0$ to be the number

$$(2.1) \qquad n_{s_0} = \mathrm{Tr}\left(\frac{1}{2\pi i}\int_{\gamma_{s_0,\varepsilon}} \mathcal{S}_X(s)^{-1}\mathcal{S}_X'(s)\,ds\right)$$

where $\gamma_{s_0,\varepsilon}$ is a closed curve of index 1 enclosing $s_0$ and no other singularity of $S_X(s)$. It follows from the analysis of [35], section 5 that the operator $\mathcal{S}_X(s)^{-1}\mathcal{S}_X'(s)$ has finite-rank residues at its singularities, and that the number $n_{s_0}$ is an integer.

## 3. Schottky Manifolds

We begin with the following standard construction of classical Schottky groups:

**Construction 3.1.** Let $\{C_j : j = 1, \ldots, 2g\}$ be a collection of circles on the Riemann sphere $S^2$, such that the $C_j$'s are disjoint and have mutually disjoint interiors. For each $i$, $i = 1, \ldots, g$, let $A_i$ be a Möbius transformation that takes $C_{2i-1}$ to $C_{2i}$, and takes the interior of $C_{2i-1}$ to the exterior of $C_{2i}$.

Let $\mathcal{D}$ denote the domain of $S^2$ exterior to all the circles $C_i$. The following lemma is then classical:

**Lemma 3.2.** *The $A_i$'s generate a free, discrete group of Möbius transformations $\Gamma$. If we denote by $\Lambda$ the limit set of $\Gamma$, then a fundamental domain for the action of $\Gamma$ on $\Omega = S^2 - \Lambda$ is given by $\mathcal{D}$.*

We remark that it is easy to check that the limit set $\Lambda$ is a Cantor set of Lebesgue measure $0$ in $S^2$.

It is not difficult to see that $S = \Omega/\Gamma$ is conformally equivalent to a surface of genus $g$. Indeed, $\Omega/\Gamma$ is obtained from $\mathcal{D}$ by gluing the circle $C_{2i-1}$ to the circle $C_{2i}$ via the map $A_i$, for $i = 1, \ldots, g$, to obtain $g$ handles. It will be convenient to observe that the $g$ gluing lines give us $g$ disjoint and homotopically distinct simple closed curves on $S$ that together do not disconnect $S$. In fact, removing the $g$ curves from $S$ leaves the connected region $\mathcal{D}$.

We may now extend the action of $\Gamma$ to hyperbolic 3-space $\mathbb{H}^3$. A fundamental domain for the action of $\Gamma$ on $\mathbb{H}^3$ may be obtained by viewing $S^2$ as the boundary of the 3-disk $\mathbb{D}^3$. For each $C_i$, we consider the hemisphere $H_i$ in $\mathbb{D}^3$ whose boundary is $C_i$. The fundamental domain for the action of $\Gamma$ on $\mathbb{H}^3$ is obtained by removing, for each $i$, the component of $\mathbb{H}^3 - H_i$ containing the interior of $C_i$.

It is easily seen from this that the quotient $M = \Gamma \setminus \mathbb{H}^3$ is a handlebody of genus $g$, whose boundary is the surface $S$. Furthermore, the $g$ curves obtained as the gluing lines above become homotopically trivial in $M$, because they may be retracted along the hemispheres $H_i$.

We remark that this construction goes through with little change if one allows the $C_i$'s to be Jordan curves rather than round circles, provided that the maps $A_i : C_{2i-1} \to C_{2i}$ remain Möbius transformations. The only change that is needed in the discussion above is in the description of the fundamental domain $F$, since there is no geometrically natural replacement for the hemispheres $H_i$.

A group $\Gamma$ of Möbius transformations constructed in this way is called a *Schottky group*, and the quotient hyperbolic manifold $M$ a *Schottky manifold*.

Our first observation is that we may invert this procedure:

**Lemma 3.3.** *Let $S$ be a Riemann surface of genus $g > 1$, and let $\{\gamma_i : i = 1, \ldots, g\}$ be a family of $g$ disjoint, homotopically distinct simple closed geodesics on $S$ such that $S - \cup_i \gamma_i$ is connected.*

*Then there are $2g$ Jordan curves $C_j$ on $S^2$ with disjoint interiors, and for each $i$, $i = 1, \ldots, g$ a Möbius transformation $A_i : C_{2i-1} \to C_{2i}$, such that if $\Gamma$ is the group generated by the $A_i$'s, then*

$$S = \Gamma \setminus \Omega(\Gamma),$$

*and the curves $\gamma_i$ are homotopically trivial in $M = \Gamma \setminus \mathbb{H}^3$.*

*Proof.* Let us make a particular choice of circles in Construction 3.1 of classical Schottky groups, to obtain a fixed Schottky group $\Gamma_0$ with corresponding surface $S_0$ and curves $\{\gamma_{0,i} : i = 1, \ldots, g\}$.

We now consider the surface $S$. Since $S - \cup_i \gamma_i$ is connected, the classification of surfaces tells us that it is diffeomorphic to the sphere $S^2$ with $2g$ holes removed. Therefore, there is a smooth map $\phi : S_0 \to S$ that takes the curves $\gamma_{0,i}$ to the curves $\gamma_i$. Since $S_0$ and $S$ are compact, $\phi$ must be quasi-conformal.

We may then use $\phi$ to pull back the conformal structure on $S$ to $S_0$. Lifting this conformal structure to $\Omega \subset S^2$ and hence to all of $S^2$ gives us a measurable conformal structure on $S^2$ invariant under $\Gamma$ and quasi-conformally equivalent to the standard conformal structure. Taking this new conformal structure on $S^2$ to the standard conformal structure conjugates $\Gamma_0$ to a new Schottky group $\Gamma$ that clearly has the desired properties. $\qquad\square$

**Definition 3.4.** Given a Riemann surface $S$, we will say that a family of $g$ disjoint, homotopically distinct simple closed curves on $S$ that together do not separate $S$ *defines a Schottky structure* on $S$. Two Schottky structures are the same if the identity map on $S$ extends to a homeomorphism of the corresponding Schottky manifolds.

It follows that two Schottky structures are the same if and only if the corresponding Schottky manifolds are isometric.

It is evident that many different families of curves may define the same Schottky structure on $S$. Indeed, suppose we are given family $\{\gamma_i\}$ of curves defining a Schottky structure, i.e., a handlebody $\Gamma \setminus \mathbb{H}^3$. A second family $\{\gamma_i'\}$ will define the same Schottky structure if and only if each of the $\gamma_i'$ is homotopically trivial in $M = \Gamma \setminus \mathbb{H}^3$.

We would like a method to determine, given two families $\{\gamma_i\}$ and $\{\gamma_i'\}$, when they define distinct Schottky structures. This is accomplished by the following:

**Lemma 3.5.** *Suppose that there exist $i$ and $j$ such that the homological intersection number $\gamma_i \cdot \gamma_j'$ is non-zero. Then $\{\gamma_i\}$ and $\{\gamma_i'\}$ determine distinct Schottky structures.*

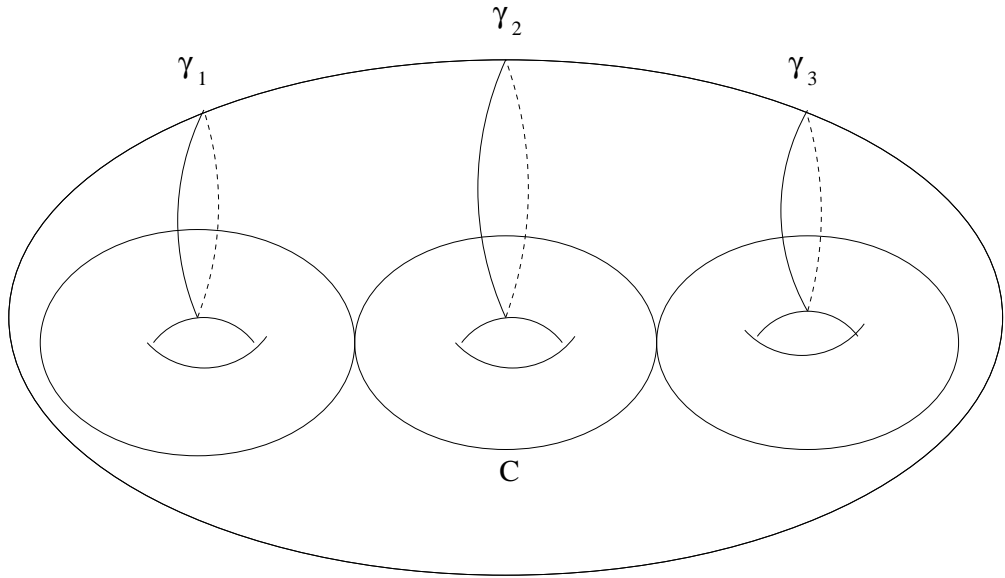

FIGURE 1. The curves $\gamma_i$

*Proof.* According to the argument in Lemma 3.3, we may assume that the $\gamma_i$'s form the standard collection of $g$ curves shown in Figure 1 above. The handlebody defined by the curves $\gamma_i$ is then the interior of $S$ viewed as a surface in $\mathbb{R}^3$. We have also drawn the core $C$ of the handlebody. This is a bouquet of $g$ circles inside the handlebody, such that the handlebody retracts onto $C$.

Let us now consider a curve $\gamma_j'$. The homological intersection number $\gamma_i \cdot \gamma_j'$ has the following interpretation *vis à vis* $C$: Let us first retract $\gamma_j'$ onto $C$, and then collapse all the circles in $C$ to points except the $i$-th circle. Then $\gamma_i \cdot \gamma_j'$ is precisely the winding number of this curve about this circle.

It follows immediately that if $\gamma_i \cdot \gamma_j' \neq 0$, then $\gamma_j'$ cannot be homotopically trivial in the handlebody. $\qquad\blacksquare$

We close this section with the following observation:

**Lemma 3.6.** *Let $S$ be a Riemann surface of genus $g$, with a family $\gamma_i$ of $g$ curves defining a Schottky structure on $S$, and let $S'$ be a $k$-fold covering surface of $S$. Then the Schottky structure on $S$ lifts to a Schottky structure on $S'$ if and only if the inverse image of each of the curves $\gamma_i$ consists of $k$ disjoint closed curves. In particular, if $S'$ is a normal covering of $S$ with covering group $G$, the Schottky structure will lift if and only if the image of each of the $\gamma_i$'s in $G$ is trivial.*

*Proof.* To say that the Schottky structure lifts is to say that there is a Schottky manifold $M'$ such that $S'$ is the boundary of $M'$, and $M'$ covers $M$.

If such an $M'$ exists, then each $\gamma_i$, being homotopically trivial in $M$, lifts to $k$ closed curves in $M'$, and hence in $S'$. This establishes the necessity of the condition.

If the inverse image of each of the $\gamma_i$'s consists of $k$ closed curves, then we obtain a family of $kg$ disjoint and homotopically distinct simple closed curves, the union of which separates $S'$ into $k$ components. Since the genus of $S'$ is $k(g-1) + 1$, it is a simple matter to delete $k - 1$ curves from this family to get a set of curves that does not disconnect $S'$.

It is easily seen that this family of curves defines a Schottky structure with the desired properties. $\qquad\blacksquare$

## 4. Transplantation on Convex Co-Compact Hyperbolic Manifolds

Given the discussion of the scattering operator in §2, it is natural to say that two manifolds are *isoscattering* if they have the same set of scattering poles and multiplicities. We want to define a stronger notion of "isoscattering" for manifolds that implies equality of the set of scattering poles together with their multiplicities.

**Definition 4.1.** Two convex co-compact hyperbolic manifolds $X_1$ and $X_2$ are called *strongly isoscattering* if for some choice of defining functions on $X_1$ and $X_2$, there is a map $\mathcal{T}_\partial : C^\infty(\partial \bar{X}_1) \to C^\infty(\partial \bar{X}_2)$ independent of $s$ that intertwines $\mathcal{S}_{X_1}(s)$ and $\mathcal{S}_{X_2}(s)$ for all $s$ in the joint domain of analyticity of $\mathcal{S}_{X_1}(s)$ and $\mathcal{S}_{X_2}(s)$.

If $X_1$ and $X_2$ are strongly isoscattering, then the set of poles of the scattering operators is the same and the multiplicities are also the same since the residues of $\mathcal{S}_{X_1}(s)^{-1}\mathcal{S}'_{X_1}(s)$ and $\mathcal{S}_{X_2}(s)^{-1}\mathcal{S}'_{X_2}(s)$ are similar finite rank operators, and hence have the same trace, at any singularity.

The main goal of this section, Theorem 4.4, is to adapt the construction of Sunada [45] to provide examples of strongly isoscattering manifolds. To that end, we need the following definition:

**Definition 4.2.** Let $(G, H_1, H_2)$ be a triple of finite groups, with $H_i$ a subgroup of $G$ for $i = 1, 2$.

$(G, H_1, H_2)$ *satisfies the Sunada condition* if, for each $g \in G$, we have

$$(4.1) \qquad \#([g] \cap H_1) = \#([g] \cap H_2)$$

where $[g]$ denotes the conjugacy class of $g$ in $G$.

The following well-known lemma gives an alternative description of the Sunada condition:

**Lemma 4.3.** *The triple* $(G, H_1, H_2)$ *satisfies the Sunada condition 4.1 if and only if the spaces* $L^2(G/H_1)$ *and* $L^2(G/H_2)$ *are equivalent as* $G$-*modules.*

The lemma can be proved readily by showing that the left $G$-representations $L^2(G/H_1)$ and $L^2(G/H_2)$ have the same characters if and only if (4.1) holds.

We now claim:

**Theorem 4.4.** *Suppose that* $(G, H_1, H_2)$ *satisfies the Sunada condition (4.1). If* $X$ *is a convex cocompact hyperbolic manifold, and* $\phi : \pi_1(X) \to G$ *an onto homomorphism, let* $X_i$ *be the covering of* $X$ *with*

$$\pi_1(X_i) = \phi^{-1}(H_i).$$

*Then* $X_1$ *and* $X_2$ *are strongly isoscattering.*

*Proof.* The idea of the proof is to adapt the transplantation method of Buser [10] and Bérard [2] to our setting. To that end, let us review this method briefly.

We begin with the observation that the group $G$ acts on the space $L^2(G)$ on both the left and the right. We may identify $L^2(G/H_i)$ with the $H_i$-invariant subspace of $L^2(G)$ under the right action of $G$, which we denote by $(L^2(G))^{H_i}$. The condition of Lemma 4.3 may then be rewritten as the existence of a vector space isomorphism

$$T : (L^2(G))^{H_1} \to (L^2(G))^{H_2}$$

that is equivariant with respect to the left action of $G$.

Let us now denote by $X_e$ the covering of $X$ with $\pi_1(X_e) = \phi^{-1}(\mathrm{id})$. We may identify functions on $X_i$ with functions on $X_e$ that are invariant under $H_i$. We observe as well that functions on $\partial \bar{X}_i$ may be similarly identified with functions on $\partial \bar{X}_e$ invariant under $H_i$.

We now pick a defining function $\rho$ on $X$, and lift it to defining functions on $X_1, X_2$, and $X_e$, that we will continue to denote by $\rho$. We may now define the scattering operators $S_i(s)$ and the Poisson operators $\mathcal{P}_i^s$ as well as $S_e(s)$ and $\mathcal{P}_e^s$ with respect to these defining functions, and we observe that with these definitions, the operators $S_i(s)$ and $\mathcal{P}_i^s$ are just the restrictions of $S_e(s)$ and $\mathcal{P}_e^s$ to the $H_i$-invariant functions.

We now pick a $G$-invariant partition of unity $\psi_g$ on $X_e$ in the following way: we pick a fundamental domain $F$ for $X$ in $X_e$, and choose $\psi_e^0$ to be a smooth function that is positive on $F$ and has support inside a small neighborhood of $F$. We may then set

$$\psi_g = \frac{g \cdot \psi_e^0}{\sum_{g' \in G} (g') \cdot \psi_e^0}.$$

For any function $f$ on $X_e$, we may decompose $f$ as

$$f = \sum_g f_g = \sum_g \psi_g \cdot f,$$

and if $f$ is $H_i$-invariant, then $f_g = f_{gh}$ for $h \in H_i$. We may now extend the map $T : (L^2(G))^{H_1} \to (L^2(G))^{H_2}$ to a map

$$\mathcal{T} : (C^\infty(X_e))^{H_1} \to (C^\infty(X_e))^{H_2},$$

and similarly for

$$\mathcal{T}_\partial : (C^\infty(\partial \bar{X}_e))^{H_1} \to (C^\infty(\partial \bar{X}_e))^{H_2}.$$

We remark that the $G$-equivariance of $T$ insures that $\mathcal{T}$ maps eigenfunctions of the Laplacian to eigenfunctions of the Laplacian.

We will show:

**Lemma 4.5.** *Let $s$ be chosen with $\Re(s) = 1$, $s \neq 1$. Then:*

(i) $\mathcal{T}\mathcal{P}_e^s = \mathcal{P}_e^s\mathcal{T}_\partial$.
(ii) $\mathcal{T}_\partial S_e(s) = S_e(s)\mathcal{T}_\partial$.

*Proof.* Indeed, let $f$ be an $H_1$-invariant function on $X_e$. Then the right-hand side in (i) applied to $f$ gives an eigenfunction of the Laplacian with eigenvalue $(s)(2-s)$ that has the asymptotics

$$S_e(s)\mathcal{T}_\partial(f)\rho^s + \mathcal{T}_\partial(f)\rho^{2-s} + O(\rho^2)$$

as $\rho \to 0$.

But the left-hand side of (i) applied to $f$ is also an eigenfunction of the Laplacian with the same eigenvalue, because $\mathcal{T}$ preserves eigenfunctions, and has the asymptotics

$$\mathcal{T}_\partial S_e(s)(f)\rho^s + \mathcal{T}_\partial(f)\rho^{2-s} + O(\rho^2).$$

Theorem 2.1 then tells us that the two sides of (i) are equal, and that

$$\mathcal{T}_\partial S_e(f) = S_e\mathcal{T}_\partial(f),$$

establishing (ii). □

We now complete the proof of Theorem 4.4 as follows: we have shown that for all $s$ with $\Re(s) = 1$, $s \neq 1$, $\mathcal{T}_\partial$ intertwines $S_1(s)$ and $S_2(s)$. The result now follows for all $s$ by analytic continuation. □

## 5. Construction of Isoscattering Manifolds

In this section, we will use Theorem 4.4 to show:

**Theorem 5.1.** *There exist $\bar{X}_1$ and $\bar{X}_2$ such that:*

(i) $X_1$ *is strongly isoscattering to $X_2$,*
(ii) $\partial\bar{X}_1$ *is conformally equivalent to $\partial\bar{X}_2$, and*
(iii) $X_1$ *is not isometric to $X_2$.*

The proof of Theorem 5.1 is closely modeled on the proof of the following result from [8]:

**Theorem 5.2** ([8])**.** *There exists a Riemann surface $S$ and functions $q_1$ and $q_2$ on $S$, such that the Schrödinger operators $\Delta_S + q_1$ and $\Delta_S + q_2$ are isospectral, but $q_1 \neq q_2$.*

We recall the idea of the proof of Theorem 5.2 briefly.

Suppose that we can find a quadruple $(G', G, K_1, K_2)$ of finite groups with the following properties:

(i) $G$ is of index 2 in $G'$.
(ii) $(G, K_1, K_2)$ satisfies the Sunada condition (4.1).
(iii) $K_1$ is conjugate to $K_2$ in $G'$.
(iv) $K_1$ is not conjugate to $K_2$ in $G$.

Now let $S_{G'}$ be a Riemann surface and $\phi : \pi_1(S_{G'}) \to G'$ an onto homomorphism. Denoting by $S_G$, $S_{K_1}$, and $S_{K_2}$ the corresponding coverings, we observe that $S_{K_1}$ is isometric to $S_{K_2}$, since $K_1$ is conjugate to $K_2$ in $G'$. Let $\zeta : S_{K_1} \to S_{K_2}$ be this isometry.

If we now choose a "sufficiently bumpy" function $q_G$ on $S_G$ that is not invariant under the action of $G'/G$ on $S_G$, and set $q_{K_1}$ and $q_{K_2}$ equal to the lifts of $q_G$ to $S_{K_1}$ and $S_{K_2}$ respectively, then the Sunada Theorem tells us that $\Delta_{S_{K_1}} + q_{K_1}$ is isospectral to $\Delta_{S_{K_2}} + q_{K_2}$.

We may now choose $S = S_{K_1}$, $q_1 = q_{K_1}$, and $q_2 = \zeta^*(g_{K_2})$ to establish the theorem.

An example of such a quadruple is (see [8] for a proof):

**Example 5.3.**    • $G = PSL(3, \mathbb{Z}_p)$ for $p$ a prime

• $K_1 = \begin{pmatrix} * & * & * \\ 0 & * & * \\ 0 & * & * \end{pmatrix}$,

• $K_2 = \begin{pmatrix} * & 0 & 0 \\ * & * & * \\ * & * & * \end{pmatrix}$,

• $G'$ is the semi-direct product of $G$ with the map $A \mapsto (A^{-1})^t$ acting on $G$.

We remark that $G$ as above can be generated by two elements, and so $G'$ can be generated by three elements.

We will now establish Theorem 5.1 by mimicking the proof of Theorem 5.2, in the following way: we begin with a quadruple $(G', G, K_1, K_2)$ as in Example 5.3, and a Riemann surface $S_{G'}$ of genus 3 as in Figure 1. We will pick a homomorphism $\pi_1(S_{G'}) \to G'$ in the following way: we first map $S_{G'}$ to the core $C$, and then map $\pi_1(C) = \mathbb{Z} * \mathbb{Z} * \mathbb{Z}$ to $G'$ so that going around the two left-hand loops correspond to generators of $G$, while going around the right-hand loop corresponds an element of $G'$ that maps to the non-trivial element, $\tau$, of $G'/G$.

Denote by $S_G$ the double cover of $S_{G'}$ corresponding to $G'/G$. We observe that $S_G$ is a surface of genus 5, and that $\tau$, the non-trivial element of $G'/G$, acts on $S_G$ by a 180° rotation about the central hole, as in Figure 2 above. We have also drawn on $S_G$ five curves $\gamma_1', \ldots, \gamma_5'$, that are lifts of the curves $\gamma_1, \gamma_2, \gamma_3$ of Figure 1 and a curve $D$ that winds once around the right-hand hole, meeting $\gamma_5'$ once, which we will use later.

Let $S_{K_1}$ and $S_{K_2}$ be the coverings of $S_G$ as constructed above. The idea is now to pick a Schottky structure on $S_G$ that will lift to Schottky structures on $S_{K_1}$ and $S_{K_2}$, but so that the isometry $\phi : S_{K_1} \to S_{K_2}$ sends the Schottky structure on $S_{K_1}$ to an inequivalent Schottky structure on $S_{K_2}$. In effect, the Schottky structure on $S_G$ will play the role of a "sufficiently bumpy" potential function. Once this is done, it follows from Theorem 4.4 that the Schottky manifolds are strongly isoscattering.

Using Lemmas 3.5 and 3.6 to define such a Schottky structure on $S_G$, we pick five curves $\{\gamma_1, \ldots, \gamma_5\}$ on $S_g$ with the following properties:

(i) The $\gamma_i$'s are disjoint simple closed curves that jointly do not separate $S_G$.

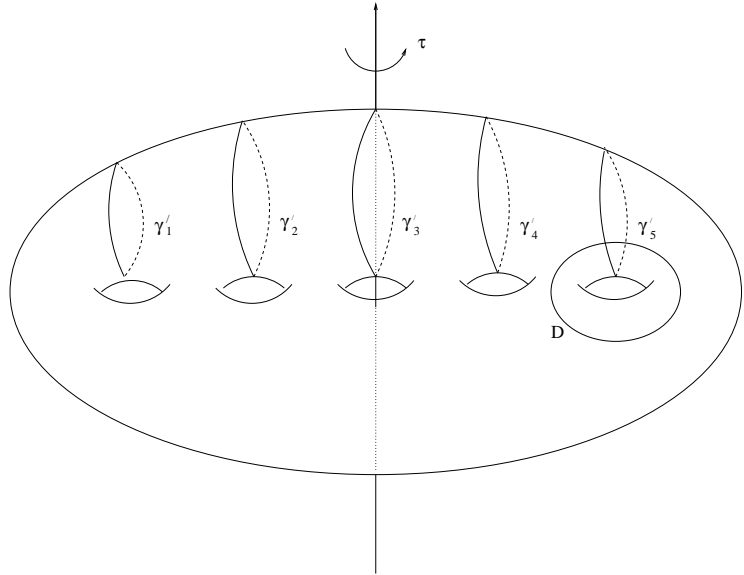

FIGURE 2. The surface $S_G$

(ii) The image of $\gamma_i$ in $G$ is trivial for all $i$.

(iii) For some $i$ and $j$, $\tau(\gamma_i)$ and $\gamma_j$ have non-trivial homological intersection number.

One such choice of curves $\gamma_1, \ldots, \gamma_5$ is shown in Figure 3 below, and may be described as follows: $\gamma_1, \ldots, \gamma_4$ are precisely the curves $\gamma_1', \ldots, \gamma_4'$ of Figure 2. Their images in $G$ are trivial because their images in the core $C$ are homotopically trivial.

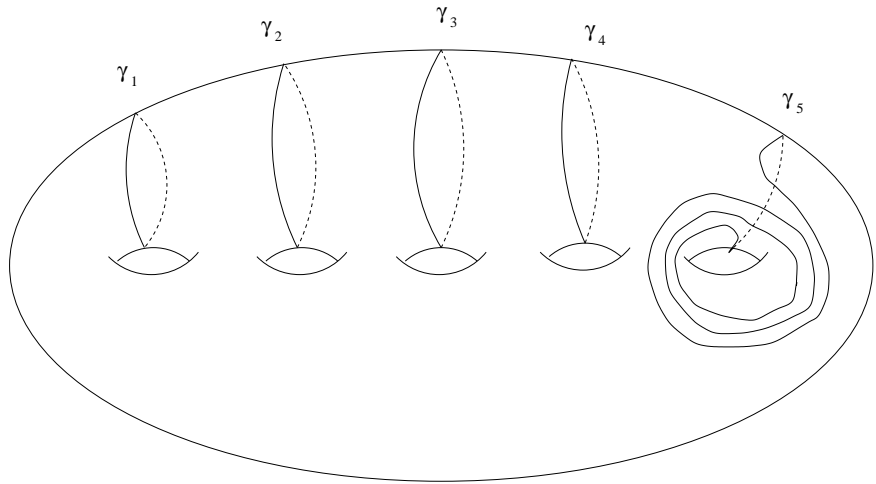

FIGURE 3. The surface $S_G$ with curves defining a distinct Schottky structure

The curve $\gamma_5$ is chosen so that it winds around the right-hand hole $k$ times, where $k$ will be chosen later, and then closes up by traveling along $\gamma_5'$. More precisely, $\gamma_5$ is the $k$-fold Dehn twist of $\gamma_5^l$ along the curve $D$. From this description, we readily see that $\{\gamma_1, \ldots, \gamma_5\}$ do not jointly separate $S_G$.

It is easily seen that the image of $\gamma_5$ in $G$ is just $(\tau(g_1))^k$, where $g_1$ is the image in $g$ of going around the left-hand loop in $C$. Since $G$ is a finite group, we may choose $k$ positive so that $(\tau(g_1))^k$ is trivial.

It is also easy to see that the homological intersection number $\tau(\gamma_1) \cdot \gamma_5$ is precisely $k$. The curves $\{\gamma_1, \ldots, \gamma_5\}$ thus have the desired properties.

In order to complete the proof of Theorem 5.1, we must show that the two Schottky manifolds $M_{K_1}$ and $M_{K_2}$ are not isometric. But any isometry of $M_{K_1}$ and $M_{K_2}$ must take $S_{K_1}$ to $S_{K_2}$ via an isometry that takes the Schottky structure on $S_{K_1}$ to the Schottky structure on $S_{K_2}$. We have constructed the Schottky structures so that the standard isometry between $S_{K_1}$ and $S_{K_2}$ does not preserve Schottky structures. We now want to choose $S_{G'}$ so that no unwanted isometries appear. But a result of Greenberg [16] and Margulis [26], as used in [45], guarantees that for a generic choice of $S_{G'}$, the commensurability group of $S_{G'}$ is trivial. It follows that any isometry from $S_{K_1}$ to $S_{K_2}$ must cover the involution $\tau : S_G \to S_G$.

This completes the proof of Theorem 5.1.

$\square$

DEPARTMENT OF MATHEMATICS, TECHNION-ISRAEL INSTITUTE OF TECHNOLOGY, HAIFA, ISRAEL
*E-mail address:* `rbrooks@techunix.technion.ac.il`

DEPARTMENT OF MATHEMATICS AND STATISTICS, TEXAS TECH UNIVERSITY, LUBBOCK, TEXAS 79409, U. S. A. CURRENT ADDRESS: DEPARTMENT OF MATHEMATICS, UNIVERSITY OF KENTUCKY, LEXINGTON, KENTUCKY, 40506–0027, U. S. A.
*E-mail address:* `gornet@math.ttu.edu`

DEPARTMENT OF MATHEMATICS, UNIVERSITY OF KENTUCKY, LEXINGTON, KENTUCKY, 40506–0027, U. S. A.
*E-mail address:* `perry@ms.uky.edu`